%
%
%
%
%
%
%
%
%
\scrollmode
\magnification=\magstep1

\hoffset=1.5cm
\hsize=12cm

\def\demo#1:{\par\medskip\noindent\it{#1}. \rm}
\def\ni{\noindent}               
\def\ll{\leftline}
\def\cl{\centerline}

\def\begin{\ll{}\vskip 10mm \nopagenumbers}  
\def\pn{\footline={\hss\tenrm\folio\hss}}   
\def\ii#1{\itemitem{#1}}

%
%
\outer\def\beginsection#1\par{\bigskip
  \message{#1}\leftline{\bf\&#1}
  \nobreak\smallskip\vskip-\parskip\noindent}

%
%
\outer\def\proclaim#1:#2\par{\medbreak\vskip-\parskip
    \noindent{\bf#1.\enspace}{\sl#2}
  \ifdim\lastskip<\medskipamount \removelastskip\penalty55\medskip\fi}

\def\endpr{\hfill $\spadesuit$ \medskip}

%
%
%
%


%
%

\def\R{{\rm I\kern-0.2em R\kern0.2em \kern-0.2em}}
\def\N{{\rm I\kern-0.2em N\kern0.2em \kern-0.2em}}
\def\P{{\rm I\kern-0.2em P\kern0.2em \kern-0.2em}}
\def\B{{\rm I\kern-0.2em B\kern0.2em \kern-0.2em}}
\def\C{{\rm C\kern-.4em {\vrule height1.4ex width.08em depth-.04ex}\;}}
\def\CP{\C\P}

%
%
%
%

\def\cJ{{\cal J}}

\def\cO{{\cal O}}

\def\cR{{\cal R}}
\def\cS{{\cal S}}

%
%
%
\def\a{\alpha}
\def\b{\beta}

\def\d{\delta}




%
%
%
%
\def\bs{\backslash}              

%
%
\def\dim{{\rm dim}\,}                    
\def\holo{holomorphic}                   
\def\nbd{neighborhood}                   

\def\ss{\subset\!\subset}                

\def\hvf{holomorphic vector field}
\def\hvb{holomorphic vector bundle}

\def\hra{\hookrightarrow}


\def\wt{\widetilde}
\def\wh{\widehat}


\begin
\cl{\bf THE OKA PRINCIPLE, LIFTING OF HOLOMORPHIC MAPS}
\medskip
\cl{\bf AND REMOVABILITY OF INTERSECTIONS}
\bigskip
\cl{Franc Forstneri\v c}
\bigskip\medskip\rm

\ni\bf Abstract. \rm
In section 1 we survey results on the Oka principle
for sections of \holo\ submersions over Stein manifolds.
In section 2 we apply these results to the lifting problem 
for holomorphic mappings, and in section 3 we apply
them to removability of intersections of holomorphic maps 
from Stein manifolds with closed complex subvarieties
of the target space.
\medskip

\beginsection 1. The Oka principle for sections of holomorphic submersions.

Our main references for this section are the papers by
Grauert [Gra], Cartan [Car], Gromov [Gro], and Prezelj
and the author [FP1, FP2, FP3]. Let $h \colon Z\to X$ be a
holomorphic map of a complex manifold
$Z$ onto a complex manifold $X$. A {\it section} of $h$ is a \holo\
map $g\colon X\to Z$ such that $hg=id_X$. We shall say
$h$ satisfies the basic {\bf Oka principle}
if the following holds:

%
%
\proclaim The basic Oka principle:
Given any continuous section $f_0\colon X\to Z$ of
$h\colon Z\to X$, there exists a homotopy of continuous
sections $f_t\colon X\to Z$ ($t\in [0,1]$)
such that $f_1$ is \holo\ on $X$.

Stronger versions concern parametrized families of sections 
with parameter in a compact Hausdorff space (see [Gro] and [FP2]).
The parametric form is essentially equivalent to the validity of

%
%
\proclaim The strong Oka principle:
The inclusion $\iota\colon \Gamma_h(X,Z)\hra \Gamma_c(X,Z)$ of
the space of holomorphic sections of $h\colon Z\to X$ to
the space of continuous sections is a weak homotopy equivalence,
that is, $\iota$ induces an isomorphism of all fundamental groups
of the two spaces (when endowed with the compact-open topology).

In 1957 H.\ Grauert established the strong Oka principle 
for certain types of holomorphic fiber bundles over Stein 
manifolds (and reduced Stein spaces) whose fiber is a complex Lie group 
or a complex homogeneous space [Gra, Car]. Grauert's result implies that  
the classification of \hvb s over Stein spaces agrees
with their topological classification (see the survey [Lei]). 
In 1986 a different proof of Grauert's theorem
was given by Henkin and Leiterer [HL1, HL2].
In 1989 M.\ Gromov [Gro] gave an important extension 
of Grauert's theorem which we now explain.

%
%
%
%
\proclaim Definition 1: {\rm [Gro, 1.1.B])}
A {\bf spray} associated to a \holo\ submersion $h\colon Z \to X$
is a triple $(E,p,s)$, where $p\colon E\to Z$ is a \hvb\ and
$s\colon E\to Z$ is a \holo\ map such that for each $z\in Z$ we
have (i) $s(E_z) \subset Z_{h(z)}$ (equivalently, $h s=h p$),
(ii) $s(0_z)=z$, and (iii) the derivative
$ds \colon T_{0_z} E \to T_z Z$ maps the
subspace $E_z \subset T_{0_z} E$ surjectively
onto $VT_z(Z):=\ker dh_z$. A spray on a complex
manifold $Z$ is by definition a spray associated to the 
trivial submersion $Z\to point$.

\demo Example:
Suppose that $V_1,V_2,\cdots,V_N$ are
$\C$-complete \hvf s  on $Z$ which are vertical
with respect to a submersion $h\colon Z\to X$
(i.e., they are tangent to the fibers of $h$).
Let $\phi^t_j$ denote the flow of $V_j$.
$\C$-completeness of $V_j$ means that $\phi_j^t$ is
defined for all complex values of the time parameter $t$
and hence $\{\phi^t_j \colon t\in \C\}$ is a complex
one-parameter group of \holo\ automorphisms of $Z$
preserving the fibers of $h$. The map
$$ s\colon Z\times \C^N \to Z, \quad
    s(z;t_1,\ldots,t_N) =
      \phi_1^{t_1}\phi_2^{t_2}
        \cdots \phi_N^{t_N}(z),
$$
satisfies properties (i) and (ii) of sprays, and it 
also satisfies (iii) when the vector fields 
$V_1,\ldots,V_N$ span the vertical tangent space 
$VT_z Z$ at each point $z\in Z$.
\endpr

\pn

The following is the Main Theorem in [Gro]
and Theorem 1.2 in [FP2].

\proclaim 1.1 Theorem: Let $h\colon Z\to X$ be a \holo\
submersion of a complex manifold $Z$ onto a Stein manifold
$X$. If each point $x\in X$ has an open \nbd\ $U\subset X$
such that the submersion $h\colon Z_U=h^{-1}(U)\to U$
admits a spray then the strong Oka principle holds for
sections of $h$. This holds in particular if $h\colon Z\to X$
is a \holo\ fiber bundle whose fiber admits a spray.

A complete proof and some further extensions can be found
in the papers [FP1, FP2, FP3]; see also the recent preprint 
[Lar] by F.\ L\'arusson. Perhaps the most important application 
(besides Grauert's classification of \hvb s) has been
the embedding theorem for Stein manifolds into euclidean 
spaces of minimal dimension, due to Eliashberg and Gromov [EGr]
and Sch\"urmann [Sch] (see also [Pre] for an
embedding--interpolation theorem).
We give another application in section 3.

\medskip\ni\bf  Additions to Theorem 1.1. \rm
\ni A1. If $K\ss X$ is a compact holomorphically convex
set and the initial section $f_0$ (or family of sections)
is \holo\ in an open set containing $K$ then the homotopy
can be chosen to approximate $f_0$ uniformly on $K$
(Theorem 1.5 in [FP2]).

\medskip
\ni A2.\ If $X_0$ is a closed complex subvariety in $X$
and the initial section $f_0\colon X\to Z$ is \holo\
on $X_0$, we can choose the homotopy $f_t$ ($t\in [0,1]$)
as above such that $f_t|_{X_0}=f_0|_{X_0}$ for all $t$
(Theorem 1.4 in [FP3]). In this case it suffices to assume 
that $h\colon Z\to X$ admits a spray over a small neighborhood 
of any point $x\in X\bs X_0$ 
(no spray is needed over points in $X_0$).

\medskip
\ni A3.\ Theorem 1.1 also holds if $X$ is a reduced Stein 
space which is stratified by a finite descending chain of 
closed complex subspaces 
$X=X_0\supset X_1\supset\cdots\supset X_m=\emptyset$
such that each stratum $S_k=X_k\bs X_{k+1}$ is non-singular and 
the restriction of the submersion $h\colon Z \to X$  
to $S_k$ admits a spray over an open \nbd\ (in $S_k$)
of any point $x\in S_k$; see the Appendix to [FP3].
\endpr

We now give a small extension of (A2), and therefore of 
Theorem 1.4 in [FP3]; we shall need this in section 3 below.

\proclaim Definition 2:
Let $h\colon Z\to X$ be  a \holo\ submersion, $X_0\subset X$
a closed complex subvariety of $X$, and $\cS\subset \cO_X$ a 
coherent sheaf of analytic ideals with support $X_0$
(i.e., $\cS_x=\cO_{X,x}$ precisely when $x\in X\bs X_0$).
We say that local \holo\ sections $f_0$ and $f_1$ of
$h\colon Z\to X$ in a \nbd\ of a point $x\in X_0$ are $\cS$-tangent
at $x$ if there is a \nbd\ $V \subset Z$ of the point
$z=f_0(x)=f_1(x)\in Z$  and a biholomorphic map
$\phi\colon V\to V'\subset \C^N$ ($N=\dim_{\C}Z$)
such that the germ at $x$ of (any component of) the map
$\phi f_0-\phi f_1\colon U\to\C^N$
belongs to ${\cS}_x$. If $f_0$ and $f_1$ are \holo\
in an open set containing $X_0$ and $\cS$-tangent 
at each $x\in X_0$, we say that $f_0$ and $f_1$
are $\cS$-tangent and write $\d(f_0,f_1)\in \cS$.

It is easily seen that the property of being
$\cS$-tangent is independent of the choice of
local coordinates on $Z$.

\proclaim 1.2 Theorem:
Let $h\colon Z\to X$ be a \holo\ submersion onto a Stein manifold
$X$, let $X_0$ be a closed complex subvariety of $X$,
and assume that $h$ admits a spray over a small \nbd\ of
any point $x\in X\bs X_0$. Let $\cS \subset\cO_X$ be a coherent 
sheaf of analytic ideals with support $X_0$.
Given a continuous section $f_0\colon X\to Z$ which is
\holo\ in a \nbd\ of $X_0$, there is a homotopy of continuous
sections $f_t\colon X\to Z$ ($t\in [0,1]$) such that for each
$t\in [0,1]$, $f_t$ is \holo\ in a \nbd\ of $X_0$ and
satisfies $\d(f_0,f_t)\in \cS$, and the section $f_1$
is \holo\ on $X$.

The analogous result holds for parametrized families of sections 
and with uniform approximation on compact holomorphically
convex subsets of $X$.  Furthermore, the result holds if $X$ is a 
reduced Stein space. Theorem 1.2 shows that the validity 
of the Oka principle for sections of a \holo\ submersion 
$h\colon Z\to X$ extends from complex subspace $X_0\subset X$ 
to all of $X$ provided that $h$ admits a spray over a 
\nbd\ of any point in $X\bs X_0$.

\demo Proof: Theorem 1.2 was proved in [FP3] 
in the case when $\cS$ is an integral power of the ideal sheaf
of $X_0$. To prove theorem 1.2 in general the following modification 
must be made. On line 4 in the proof of Theorem 5.2 in [FP3] 
we chose finitely many \holo\ functions $h_1,\ldots, h_m$ on $X$ 
such that
$$ X_0=\{x\in X\colon h_j(x)=0,\ 1\le j\le m\} \eqno(1.1) $$
and each $h_j$ vanishes to a given order $r \in \N$
on $X_0$. We now replace this by the requirement
that {\it the functions $h_j$ be sections of the sheaf $\cS$}.
(A finite collection of sections of $\cS$ satisfying (1.1)
can be constructed by using Cartan's Theorem A; one
inductively lowers the dimension of the superfluous
irreducible components of their common zero set.)
Once this replacement is made, the proof of Theorem 1.4
in [FP3, section 6] gives theorem 1.2 above, and no other
changes are needed.
\endpr

We give an important special case of theorem 1.2.
Let $h\colon E\to X$ be a \hvb\ of rank $q$ over a
reduced Stein space $X$. For each $x\in X$ we denote by
$\wh E_x \approx \CP^q$ the compactification of
the fiber $E_x \cong \C^q$ obtained by adding to $E_x$
the hyperplane at infinity $\Lambda_x\approx\CP^{q-1}$.
The resulting fiber bundle $\wh h\colon \wh E \to X$
with fibers $\wh E_x\approx \CP^q$ is again \holo.
(The essential observation is that the transition maps,
which are $\C$-linear automorphisms of fibers $E_x$, 
extend to projective linear automorphisms of $\hat E_x$).

\proclaim 1.3 Corollary:
{\bf (Avoiding subvarieties by \holo\ sections.)}
Let $h\colon E\to X$ be a \hvb\ of rank $q$ over a reduced
Stein space $X$ and let $\wh h\colon \wh E\to X$
be the associated bundle with fiber $\CP^q$.
Let $\wh \Sigma\subset \wh E$ be a closed complex subvariety
and set $\Sigma =\wh \Sigma\cap E$.
If for each $x\in X$ the fiber $\Sigma_x=\Sigma\cap E_x$
is of complex codimension at least two in $E_x$ then
the strong Oka principle holds for sections of
$h\colon E\bs \Sigma\to X$.

Notice that the fibers $\Sigma_x$ are algebraic by
Chow's theorem. Under stronger hypothesis on $\Sigma$ 
this is Corollary 1.8 in [FP2].

\demo Proof of corollary 1.3:
In [FP2, Lemma 8.3] the following was shown:
{\it If $\Lambda_x \not\subset \wh \Sigma_x$ for some $x\in X$ 
then there is an open \nbd\ $U\subset X$ of $x$ such that the 
submersion $h\colon h^{-1}(U) \bs \Sigma \to U$ admits
a spray.} Granted this lemma we prove the corollary by
induction on $n=\dim_{\C} X$ as follows. For $n=0$ the result 
is trivial. Assume now that $n\ge 1$ and that the result holds 
over Stein spaces of dimension $<n$. Let $\Lambda\subset \wh E$ 
be the complex hypersurface with fibers $\Lambda_x=\wh E_x\bs E_x$.
If an irreducible component of $\Lambda$ is contained in
$\wh \Sigma$, we can remove this component from $\wh \Sigma$
without changing the assumption on $\Sigma=\wh \Sigma\cap E$.
Thus we may assume that $\wh \Sigma$ and $\Lambda$ have no
common irreducible components. The set
$$ 
	X_1	= \{x\in X\colon \Lambda_x\subset \wh \Sigma_x\} 
	\cup X_{\rm sing}
$$
is then a closed complex subspace of $X$ with $\dim X_1\le n-1$.
The restriction of the submersion $h\colon E\bs \Sigma\to X$ to $X_1$
satisfies the hypothesis of  corollary 1.3 and hence the 
conclusion holds over $X_1$. Any \holo\ section over $X_1$ (or 
any compact family of sections) extends holomorphically to an
open \nbd\ $V\subset X$ of $X_1$ such that they still avoid $\Sigma$
over $V$. If $x\in X\bs X_1$ then $\Lambda_x \not\subset \wh \Sigma_x$
and hence  the submersion $h\colon E\bs \Sigma\to X$ admits a spray 
over a \nbd\ of $x$ (see the lemma stated at the beginning of the proof). 
By theorem 1.2 the Oka principle extends from $X_1$ to $X$.
\endpr

\beginsection 2. Lifting of holomorphic mappings.

In this section we show that the results discussed in section 1
give the Oka principle for liftings of holomorphic maps. 
I wish to thank Takeo Ohsawa who asked
this question during my talk at the Hayama conference
on December 16, 2000. 

Let $h \colon Z\to X$ be a holomorphic map. 
Given a complex manifold $Y$ and a \holo\ map $f\colon Y\to X$, 
the problem is to find a \holo\ map $g\colon X\to Z$ such that 
$h g=f$; any such $g$ will be called a {\bf holomorphic lifting} 
of $f$ (with respect to $h$).

$$ \matrix{ & & Z & \cr
            & {g \atop\nearrow} & \downarrow & \!\!\!\!\!\!\!\! \! h \cr
            Y & {f\atop \longrightarrow} & X & \cr}
$$

An obvious necessary condition is the existence of a
{\bf continuous lifting}, and the main question is when
is this condition also sufficient.

%
%
\proclaim 2.1 Theorem: {\bf (The Oka principle for liftings.)}
Let $h\colon Z\to X$ be a \holo\ map.
Suppose that $Y$ is a Stein manifold, $f\colon Y\to X$
is a holomorphic map and $g_0\colon Y \to Z$ is a continuous
map such that $h  g_0=f$. Assume that for each $y\in Y$
the point $f(y)\in X$ has an open \nbd\ $U\subset X$ such
that $h\colon Z_U = h^{-1}(U) \to U$ is a \holo\ submersion
onto $U$ which admits a spray (definition 1).
Then there exists a homotopy of continuous maps 
$g_t\colon Y\to Z$ such that $h g_t=f$ for each $t\in [0,1]$
and the map $g_1$ is \holo\ on $Y$. If $g_0$ is
holomorphic on a closed complex subvariety $Y_0\subset Y$
then we can choose the homotopy $g_t$ to be fixed on $Y_0$.
If $g_0$ is \holo\ in a \nbd\ of a compact
holomorphically convex subset $K\ss Y$,
the homotopy $g_t$ can be chosen
to approximate $g_0$ uniformly on $K$.

\demo  Proof:
The main idea is to apply theorem 1.2 above
to a pull-back submersion $\wt h \colon \wt Z = f^* Z \to Y$
which is constructed such that liftings of $f$
correspond to section of $\wt h$ and sprays for $h$
induce sprays for $\wt h$. Set
$$
   \eqalignno{ \wt Z &= \{ (y,z)\colon y\in Y,\ z\in Z,\ f(y)=h(z)\}, \cr
   & \wt h(y,z)=y \in Y,\quad \sigma(y,z)=z\in Z. & (2.1) \cr}
$$
It is easily seen that $\wt Z$ is a closed complex
submanifold of $Y\times Z$, the maps $\wt h\colon \wt Z \to Y$
and $\sigma\colon\wt Z\to Z$ are \holo\ and
satisfy $f\wt h=h\sigma$.
Furthermore, since $h$ is a submersion in a \nbd\ of any
point $f(y)\in X$ for $y\in Y$, it follows that $\wt h$
is a \holo\ submersion of $\wt Z$ onto $Y$.
For any section $\wt g\colon Y\to\wt Z$ of
$\wt h \colon \wt Z\to Y$ the map
$g= \sigma \wt g\colon Y\to Z$ is an $h$-lifting
of $f$ since
$$
    hg=h(\sigma \wt g)=(h \sigma)\wt g
    =(f\wt h)\wt g = f(\wt h \wt g)= f.
$$
Moreover, we claim that any lifting $g$ of $f$ is of this form.
To see this, observe that $h(g(y))=f(y)$ for $y\in Y$ implies that the
point $\wt g(y):=(y,g(y)) \in Y\times Z$ belongs to the
subset $\wt Z \subset Y\times Z$ defined by (2.1) and
hence $\wt g\colon Y\to \wt Z$ is a section of
$\wt h \colon \wt Z\to Y$. Furthermore we have
$\sigma(\wt g(y))=\sigma(y,g(y))=g(y)$ and hence the
lifting $g$ of $f$ is indeed obtained from the section
$\wt g$. Therefore theorem 2.1 will follow from theorem 1.2
in section 1 once we prove the following.

%
%
\proclaim 2.2 Lemma: {\bf (Pulling back sprays.)}
Let $f\colon Y\to X$ and $h\colon Z\to X$ be \holo\ maps.
Assume that $U\subset X$ is an open set
such that map $h\colon Z_U= h^{-1}(U) \to U$ is a
submersion onto $U$ which admits a spray.
Then the map $\wt h\colon \wt Z=f^* Z\to Y$
defined by (2.1) is a submersion with spray over
the open set $V=f^{-1}(U)\subset Y$.

\demo Proof:
Let $(E,p,s)$ be a spray associated to the submersion
$h\colon Z_U\to U$. Set $V=f^{-1}(U) \subset Y$
and observe that $\sigma$ maps the set
$\wt Z_V=\wt h^{-1}(V)$ to $Z_U$.
Let $\wt p\colon \wt E\to \wt Z_V$ denote the pull-back
of the \hvb\ $p\colon E\to Z_U$ by the map
$\sigma \colon \wt Z_V \to Z_U$. Explicitly we have
$$
   \eqalign{ \wt E &= \{ (\wt z,e)\colon
    \wt z\in \wt Z_V,\ e\in E;\ \sigma(\wt z)= p(e)\} \cr
    &= \{(y,z,e)\colon y\in V,\ z\in Z,\ e\in E;\
       f(y)=h(z),\ p(e)=z\} \cr}
$$
and $\wt p(\wt z,e)=\wt z$. Consider the map
$\wt s\colon \wt E\to \wt Z_V$ defined by $s(y,z,e)=(y,s(e))$.
We claim that $(\wt E,\wt p,\wt s)$ is a spray associated to 
the submersion $\wt h\colon \wt Z_V\to V$.
We first check that the the map $\wt s$ is well defined.
If $(y,z,e)\in \wt E$ then $p(e)=z$ and $h(z)=f(y)$.
Since $s$ is a spray for $h$, we have $h(s(e))= h(z)=f(y)$
which shows that the point $\wt s(y,z,e)=(y,s(e)) \in Y\times Z$
belongs to the fiber $\wt Z_y$. This verifies property
(i) in definition 1. Clearly
$\wt s(y,z,0_{(y,z)})=(y,s(0_z))=(y,z)$ which verifies
property (ii) in definition 1. It is also immediate
that $\wt s$ satisfies property (iii) provided that $s$
does since the vertical derivatives of the two maps
coincide under the identifications
$\wt Z_{y} \approx Z_{f(y)}$ and $\wt E_{(y,z)} \approx E_{z}$.
\endpr

\demo Proof of theorem 2.1 for submersions with stratified sprays:
Assume that $X$ is stratified by a finite descending chain of
closed complex subvarieties
$X=X_0\supset X_1\supset\cdots\supset X_m=\emptyset$
such that for each $k=0,\ldots, m-1$ the stratum $S_k=X_k\bs X_{k+1}$
is non-singular and the submersion $h\colon Z_k=h^{-1}(S_k)\to S_k$
admits a spray over an open \nbd\ of any point $x\in S_k$.
Let $f\colon Y\to X$ be a holomorphic map.
Set $Y'_k=f^{-1}(X_k) \subset Y$. Then
$Y=Y'_0\supset Y'_1\supset\cdots\supset Y'_m=\emptyset$
is a stratification of $Y$. There are two problems:
(1) the strata $\Sigma'_k=Y'_k\bs Y'_{k+1}$ may have
singularities, and (2) the subvariety $Y' \subset Y$
(over which the initial map $g_0 \colon Y\to Z$
is \holo) need not be included in the sets $Y'_k$.
To rectify this we take $\wt Y_k=Y'_k\cup Y'$ to get
a stratification
$Y=\wt Y_0\supset \wt Y_1\supset\cdots\supset \wt Y_m= Y'$.
We delete any possible repetitions and pass to a refinement 
$Y=Y_0\supset Y_1\supset \cdots\supset Y_l=Y'$
in which the strata $\Sigma_k=Y_k\bs Y_{k+1}$ are
nonempty and regular.

By the construction $f$ maps each stratum
$\Sigma_k=Y_k\bs Y_{k+1}$ to a stratum 
$S_j=X_j\bs X_{j+1}$ for some $j=j(k)$. Lemma 2.2 now provides 
local stratified sprays for the pull-back submersion
$\wt h\colon \wt Z\to Y$ defined by (2.1).
More precisely, for each $y\in \Sigma_k$ we get an open
\nbd\ $V$ of $y$ in $\Sigma_k$ and a spray on the
submersion $\wt h\colon \wt h^{-1}(V) \to V$ by pulling back
a spray from a \nbd\ $U\subset S_j$ of the point
$f(y)\in S_j$ as in lemma 2.2. Hence theorem 2.1 follows
from the version of theorem 1.1 for submersion with 
stratified sprays. 
\endpr

%
%
%
\beginsection 3. Removability of intersections.

The main reference for this section is [Fo] and [FP3].
Let $Z$ be a complex manifold and $\Sigma\subset Z$
a closed complex subvariety of $Z$. Given
a complex manifold $X$ and a holomorphic map
$f\colon X\to Z$, we write
$$
    f^{-1}(\Sigma)=\{x\in X\colon f(x)\in \Sigma \}
    =Y\cup\wt Y                            \eqno(3.1)
$$
where each of the sets $Y$ and $\wt Y$ is a union of
connected components of the preimage $f^{-1}(\Sigma)$
and $Y\cap \wt Y=\emptyset$. We say that the set $\wt Y$ 
is {\bf holomorphically removable} from $f^{-1}(\Sigma)$ 
if the following holds:

%
%
\medskip\ni\bf (H-rem) \sl There is a \holo\ homotopy
$f_t\colon X\to Z$ ($t\in [0,1]$), with $f=f_0$, such that 
for each $t\in [0,1]$ the set $Y$ is a union of connected 
components of $f_t^{-1}(\Sigma)$ and we have $f_1^{-1}(\Sigma)=Y$.
\rm

\medskip
The following is clearly a necessary condition
for the validity of (H-rem):

\medskip\ni\bf (C-rem) \sl There are an open set
$U\subset X$ containing $Y$ and a homotopy of continuous
maps $\wt f_t\colon X\to Z$ ($t\in [0,1]$) with $\wt f_0=f$
such that for each $t\in [0,1]$ we have
$\wt f_t|_U=f_t|_U$, and $\wt f_1^{-1}(\Sigma)=Y$.
\rm \medskip

The validity of the Oka principle means that 
(C-rem)$\Rightarrow$(H-rem).

\proclaim 3.1 Theorem:
{\bf (The Oka principle for removability of intersections)}
Let $X$ be a Stein manifold and $f\colon X\to Z$ a \holo\
map to a complex manifold $Z$. Let $\Sigma$ be a closed complex 
subvariety of $Z$ and write $f^{-1}(\Sigma)=Y\cup\wt Y$ as in (3.1). 
Suppose that (C-rem) holds for $\wt Y$. Then:
\item{(a)} If the manifold $Z\bs \Sigma$ admits a spray,
there is a continuous homotopy $f'_t\colon X\to Z$
($t\in [0,1]$) satisfying (C-rem) such that
the map $f'_1$ is \holo\ on $X$.
\item{(b)} If both $Z$ and $Z\bs \Sigma$ admit a spray
then (C-rem)$\Rightarrow$(H-rem).

We postpone the proof for a moment and give a couple
of examples. As our first example we let $Z=\C^d$ for 
some $d\ge 1$. The map $s\colon \C^d\times \C^d\to \C^d$,
$s(z,w)=z+w$, is clearly a spray on $\C^d$.
A closed complex subvariety $\Sigma\subset \C^d$ is said
to be {\bf tame} if there is a \holo\ automorphism
$\Phi$ of $\C^d$ such that
$\Phi(\Sigma)\subset \{(z',z_d)\in \C^d \colon |z_d|\le 1+|z'|\}.
$
In particular, any algebraic subvariety of $\C^d$ is tame.
For discrete subsets of $\C^d$ this notion of tameness
coincides with the one introduced by Rosay and Rudin [RR].
The following is a corollary of theorem 3.1.

\medskip\ni\bf  3.2 Corollary. \sl
{\rm  (Theorem 3.1 in [Fo])} 
Let $\Sigma$ be a closed complex subvariety of $\C^d$
satisfying one of the following conditions:
\item{(a)} $\Sigma$ is tame and $\dim \Sigma\le d-2$;
\item{(b)} a complex Lie group acts holomorphically
and transitively on $\C^d\bs \Sigma$.

\ni Then  (C-rem)$\Rightarrow$(H-rem) holds for 
any data $(X,f,\wt Y)$ as in theorem 3.1.
\medskip\rm

\demo Proof: In each of the two cases the manifold 
$\C^d\bs \Sigma$ admits a spray.
\endpr

Corollary 3.2 applies in particular when $\Sigma=\{0\} \subset \C^d$
and gives the result of Forster and Ramspott [FR] from 1966
to the effect that the Oka principle holds in the problem 
of {\it complete intersections} on Stein manifolds.
\medskip

We mention that for each $d\ge 1$ there exist discrete sets 
$\Sigma\subset\C^d$ such that corollary 3.2 fails. In fact, 
there exists a discrete set $B\subset \C^d$ which is {\it unavoidable}, 
in the sense that every entire map $G\colon \C^d\to\C^d\bs B$ whose 
image avoids $B$ has rank $<d$ at each point (Rosay and Rudin [RR]).
Choose a point $p\in \C^d\bs B$ and set $\Sigma=B\cup\{p\}$.
Take $X=\C^d$ and let $f\colon\C^d\to\C^d$ be the identity
map $f(z)=z$. Also take $Y=\{p\}$ and $\wt Y=B$. Then clearly
(C-rem) holds but (H-rem) fails for the pair $(f,\wt Y)$
(since rank$G<d$ for a holomorphic map
$G\colon \C^d\to \C^d\bs B$ implies that $G^{-1}(p)$
contains no isolated points, and hence $p$ cannot be a
connected component of $G^{-1}(p)$).
\endpr

As our second example we let $\Sigma$ be an affine
complex subspace of $Z=\CP^d$. Since $\CP^d$ and $\CP^d\bs \Sigma$ 
are complex homogeneous spaces (the group of all affine 
complex automorphisms acts transitively on $\CP^d$, and
the group of those affine automorphisms which fix
$\Sigma$ acts transitively on $\CP^d\bs \Sigma$),
we get

\proclaim 3.3 Corollary: For any Stein manifold $X$ the
implication (C-rem)$\Rightarrow$(H-rem) holds for intersections of
holomorphic maps $X\to\CP^d$ with affine complex
subspaces of $\CP^d$.

\demo Problem: Does corollary 3.3 hold if $\Sigma\subset \CP^d$
is an algebraic subvariety of codimension at least two?
Does $\CP^d\bs \Sigma$ admit a spray for every such $\Sigma$?
In particular, does $\CP^d\bs \{p,q\}$ admit a spray?

%
%
%
%
\demo Proof of theorem 3.1 (a): By assumption there is 
a continuous homotopy $\wt f_t\colon X\to Z$
($t\in [0,1]$) satisfying (C-rem) for the given
initial map $f=\wt f_0$ and the subset
$\wt Y\subset f^{-1}(\Sigma)$. Thus $\wt f_1 \colon X\to Z$
is a continuous map which is \holo\ near the set
$Y=\wt f_1^{-1}(\Sigma)$. Assuming that $Z\bs \Sigma$ admits a 
spray, we must show that one can deform $\wt f_1$ to a 
holomorphic map $\wt f_2\colon X\to Z$ by a continuous 
homotopy $\wt f_t\colon X\to Z$ ($1\le t\le 2$)
which is \holo\ in a \nbd\ of $Y$ and satisfies
$\wt f_t^{-1}(\Sigma)=Y$ for all $t\in [1,2]$.
The homotopy $f'_t=\wt f_{2t}$ ($t\in [0,1]$)
will then satisfy part (a) in theorem 3.1.

A complete proof was given in [Fo] for $Z=\C^d$, but
unfortunately it uses the linear structure on $\C^d$.
A small modification is needed in the general case.
To simplify the notation we replace $f$ by $\tilde f_1$;
hence $f$ is \holo\ near $Y =f^{-1}(\Sigma)$.
We  define a coherent sheaf of ideals $\cR$ on $X$
which measures the order of contact of $f$ with $\Sigma$
along $Y$. Fix a point $x\in Y$ and let $z=f(x)\in \Sigma$.
Let $g_1,\ldots,g_k$ be \holo\ functions which
generate the sheaf of ideals of $\Sigma$ in some \nbd\ of
$z$ in $Z$. We take the functions $g_j \circ f$,
$1\le j\le k$, as the local generators of the sheaf
$\cR$ near $x$. Furthermore we take
$\cS=\cR\cdotp \cJ_Y^r$, where $\cJ_Y$ is the
sheaf of ideals of the subvariety $Y$ and $r$ is
a fixed positive integer. The purpose of introducing the 
sheaf $\cS$ is explained by the following lemma.

\proclaim 3.4 Lemma: {\rm (Notation as above)}
If $U\subset X$ is an open set containing $Y$
and $g\colon U\to Z$ is a \holo\ map which satisfies
$\d(f,g)\in \cS$ then there is an open set $V$,
with $Y\subset V\subset U$, such that
$\{x\in V\colon g(x)\in \Sigma\} =Y$.
The analogous conclusion holds for continuous 
families of sections with parameter in a compact space.

Lemma 3.4 is proved in [Fo, Section 3].
We continue with the proof of theorem 3.1.
We identify maps $X\to Z$ with sections of the
trivial submersion $h\colon \wt Z=X\times Z\to X$
without changing the notation. Let
$\wt Z'=\wt Z\bs \wt\Sigma=X\times (Z\bs \Sigma)$ 
where $\wt \Sigma= (X\times \Sigma)$. This this is a 
trivial submersion with a spray (since there is 
a spray on the fiber $Z\bs \Sigma$). Then $f$ is a 
section of $\wt Z$ whose image outside the subvariety 
$Y$ belongs to $\wt Z'$.

We now apply the proof of Theorem 1.4 in [FP3];
more precisely, we apply the second version
in which the sections are holomorphic near $Y$ and
the patching of sections takes place over small sets
in $X\bs Y$. There is only one apparent difficulty:
The sections have values in $\wt Z$ over points in $Y$ 
and they must have values in the smaller set $\wt Z'$ 
over points in $X\bs Y$. Fortunately this change of 
codomain is only a virtual difficulty which can be 
avoided as follows.

At a typical step of the procedure in [FP3] we have
a compact, holomorphically convex set $K\subset X$
and a pair of sections $(a,b)$ satisfying:
\item{(i)} $a\colon X\to\wt Z$ is a continuous section
which is \holo\ in an open set $\wt A\subset X$
with $Y\cup K\subset \wt A$,
\item{(ii)}  $a(x)\in \wt \Sigma $ precisely when $x\in Y$,
\item{(iii)} $b\colon \wt B\to \wt Z'$ is a \holo\ section
with values in $\wt Z'$, defined in an open \nbd\
$\wt B$ of a (small) compact set $B\ss X\bs Y$,
\item{(iv)} there is a homotopy of \holo\ sections
of $\wt Z'$ over the intersection $\wt C=\wt A\cap \wt B$
which connects $a$ to $b$.

\medskip\ni
The goal is to patch the two sections $(a,b)$ into a
single section $\wt a\colon X\to \wt Z$ which
satisfies conditions (i) and (ii) above over an open
\nbd\ of $Y\cup K\cup B$ and which approximates
$a$ uniformly on $K$. We need to assume that $(K,B)$ is 
a Cartan pair (see [FP2] and [FP3]). The patching of $a$ 
and $b$ is achieved by performing the following steps. 
We shall follow the notation in [FP3] as much as possible, 
with $X_0 =Y$. 

In the proof of Theorem 5.2 in [FP3] we constructed a 
\holo\ map $s_1 \colon V\to \wt Z$ from an open 
set $V\subset \wt A\times \C^N$  containing
$\wt A\times \{0\}$ (for some large $N\in \N$)
such that $s_1(x,0)=a(x)$ for all $x\in \wt A$ and 
the map $\xi\to s_1(x,\xi) \in \wt Z_x$ has derivative
of maximal rank at $\xi =0$ for $x\in \wt A\bs Y$.
However, for $x\in Y$ we have $s_1(x,\xi)=a(x)$ for all 
$\xi \in \C^N$ near $\xi=0$.

In the construction of $s_1$ in [FP3] we used
certain \holo\ functions $h_1,\ldots,h_m$ on $X$ which
vanish to a given order $r$ on $Y=X_0$ and whose common
zero set is precisely $Y$. The only change in our
current situation is to choose these functions $h_j$ to
be sections of the sheaf $\cS$ constructed above.
Furthermore, we let $s_2 \colon \wt B\times \C^N\to \wt Z'$
be the holomorphic map as in the proof of Theorem 5.2 in 
[FP3] (which is obtained from the spray on the submersion 
$\wt Z' \to X$ over the open set $\wt B$).

The patching of $a$ and $b$ now proceeds in two steps.
In the first step we use the holomorphic homotopy between
the two sections over $\wt C$ in order to replace
$b$ by a section of $\wt Z'$ (over a smaller \nbd\ of
$B$) which approximates $a$ sufficiently well over
a \nbd\ of $K\cap B$. This is done exactly as in [FP3].
In the second step the two sections (which are now close
over a \nbd\ of $K\cap B$) are patched by forming a transition
map $\psi$ satisfying $s_1(x,\xi)=s_2(x,\psi(x,\xi))$
(see (5.2) in [FP3]) and then solving the equation
$\psi(x,\a(x))=\b(x)$ for $x$ near $K\cap B$
(see (4.1) in [FP3]). The solutions $\a$
(which is \holo\ over a \nbd\ of $Y\cup K$)
and $\b$ (which is \holo\ over a \nbd\ of
$B$) then give a single section $\wt a$ of $\wt Z$
over an open  \nbd\ of $Y\cup K\cup B$,
defined by $\wt a(x)=s_1(x,\a(x))=s_2(x,\b(x))$.
By construction the two expressions agree for
$x\in K\cap B$.

We claim that this resulting section $\wt a$ satisfies
all requirements, provided that it approximates $a$
sufficiently well in a \nbd\ of $K$ (as we may assume
to be the case). By construction of the map $s_1$
we have $\d(\wt a,a)\in \cS$, that is, the two sections
are $\cS$-tangent along $Y$. Lemma 3.4 implies that
$\wt a(x)\in \wt Z'$ if $x$ is sufficiently
close to $Y$ but not in $Y$. If $\wt a$ approximates $a$ 
sufficiently well on $K$ it follows that $\wt a(x)\in \wt Z'$ 
for all $x\in K\bs Y$. For $x\in B$ the  same is true 
since $s_2$ has range in $\wt Z'$ and we have 
$\wt a(x) =s_2(x,\b(x))$). Finally we extend $\wt a$ 
to a continuous section over $X$ by patching it with
$a$ outside a suitable \nbd\ of $Y\cup K\cup B$.

This completes the induction step. The final section 
$f'_1$ which is \holo\ on $X$ is obtained as a locally
uniform limit of sections obtained by this procedure.
The same construction can be done for parametrized
families of sections and we get a required homotopy
$f'_t$ ($t\in [0,1]$) satisfying theorem 3.1 (a).

\demo Proof of theorem 1.3 (b):
Suppose now that $Z$ also admits a spray.
Let $\{f'_t\}$ ($t\in [0,1]$) be a homotopy satisfying 
theorem 3.1 (a). Applying theorem 1.2 in section 1 
to $\{f'_t\}$, with the sheaf $\cS$ defined at the beginning
of the proof, we obtain a two-parameter homotopy
of maps $h_{t,s}\colon X\to Z$ ($t,s\in [0,1]$) 
which are \holo\ in a \nbd\ of $Y$ and satisfy 
the following properties:
\medskip

\item{(i)} $h_{t,0}=f'_t$ for all $t\in [0,1]$,
\smallskip
\item{(ii)} $h_{0,s}=f'_0$ and $h_{1,s}=f'_1$ for all
$s\in [0,1]$ 
\smallskip
\item{(iii)} $\d(h_{t,0},h_{t,s}) \in \cS$ for all
$s,t\in [0,1]$, and
\smallskip
\item{(iv)} the map $f_t:=h_{t,1}$ is \holo\ on $X$
for each $t\in [0,1]$.

\medskip
It follows from (iii) and lemma 3.4 above that
there is a \nbd\ $V\subset X$ of $Y$ such that 
$h_{t,s}^{-1}(\Sigma)\cap V=Y$ for all $s,t\in [0,1]$. 
The homotopy $\{f_t\colon t\in [0,1]\}$, defined by
(iv) above, then satisfies (H-rem) for the initial map
$f=f_0$ and the set $\wt Y \subset f^{-1}(\Sigma)$.
\endpr

\vfill\eject
%
%
%
%
\medskip\ni\bf References. \rm
\medskip

\ii{[Car]} H.\ Cartan: Espaces fibr\'es analytiques.
Symposium Internat.\ de topo\-logia algebraica, Mexico, 97--121 (1958).
(Also in Oeuvres, vol.\ 2, Springer, New York, 1979.)

\ii{[FR]} O.\ Forster and K.\ J.\ Ramspott: Analytische
Modulgarben und Endromisb\"undel. Invent.\ Math.\ {\bf 2},
145--170 (1966).

\ii{[Fo]} F.\ Forstneri\v c: On complete intersections.
Annales Inst.\ Fourier {\bf 51} (2001).

\ii{[FP1]} F.\ Forstneri\v c and J.\ Prezelj:
Oka's principle for holomorphic fiber bundles with sprays.
Math.\ Ann.\ {\bf 317} (2000), 117-154.

\ii{[FP2]} F.\ Forstneri\v c and J.\ Prezelj:
Oka's principle for holomorphic submersions with sprays.
Preprint, 1999.

\ii{FP3]} F.\ Forstneri\v c and J.\ Prezelj:
Extending holomorphic sections from complex subvarieties.
Math.\ Z.\ {\bf 236} (2001), 43--68.

\ii{[Gra]} H.\ Grauert:
Holomorphe Funktionen mit Werten in komplexen Lieschen Gruppen.
Math.\ Ann.\ {\bf 133}, 450--472 (1957).

\ii{[Gro]} M.\ Gromov:
Oka's principle for holomorphic sections of elliptic bundles.
J.\ Amer.\ Math.\ Soc.\ {\bf 2}, 851-897 (1989).

\ii{[GR]} C.\ Gunning, H.\ Rossi:
Analytic functions of several complex variables.
Prentice--Hall, Englewood Cliffs, 1965.

\ii{[HL1}   G.\ Henkin, J.\ Leiterer:
Proof of Oka-Grauert principle without the induction over
basis dimension.
Preprint, Karl Weierstrass Institut f\"ur Mathematik,
Berlin, 1986.

\ii{[HL2]} G.\ Henkin, J.\ Leiterer:
The Oka-Grauert principle without induction over the basis dimension.
Math.\ Ann.\ {\bf 311}, 71--93 (1998).

\ii{[H\"or]} L.\ H\"ormander:
An Introduction to Complex Analysis in Several Variables, 3rd ed.
North Holland, Amsterdam, 1990.

\ii{[Lar]} F.\ L\'arusson: Excision for simplicial sheaves on the
Stein site and Gromov's Oka principle.
Preprint, 2000.

\ii{[Lei]}  J.\ Leiterer:
Holomorphic Vector Bundles and the Oka-Grauert Priciple.
Encyclopedia of Mathematical Sciences, vol. 10, 63--103;
Several Complex Variables IV, Springer,  1989.

\ii{[Pre]} J.\ Prezelj: Interpolation of embeddings of Stein
manifolds on discrete sets. Pre\-print, 1999.

\ii{[RRu]} J.-P.\ Rosay and W.\ Rudin:
Holomorphic maps from $\C^n$ to $\C^n$.
Trans.\ Amer.\ Math.\ Soc.\ {\bf 310}, 47--86 (1988)

\ii{[Sch]} J.\ Sch\"urmann:
Embeddings of Stein spaces into affine spaces of minimal dimension.
Math.\ Ann.\ {\bf 307}, 381--399 (1997).

%
%
%
%
\bigskip\ni {\it Address:} Institute of Mathematics, Physics and
Mechanics, University of Ljubljana, Jadranska 19, 1000 Ljubljana,
Slovenia

\bye